\documentclass[12pt,draft]{amsart}
\usepackage{amsmath,amsthm,latexsym,amscd,amsbsy,amssymb,fleqn,leqno}
\chardef\bslash=`\\ 

\makeatletter
\def\verbatim{\interlinepenalty\@M \@verbatim
  \leftskip\@totalleftmargin\advance\leftskip2pc
  \frenchspacing\@vobeyspaces \@xverbatim}
\makeatother
\makeatletter
  \def\dgt@k{\dg@DX=-3 \dg@DY=2 \dg@SIZE=3} 
\makeatother

\makeatletter
  \def\dgt@kk{\dg@DX=3 \dg@DY=-1 \dg@SIZE=3}%
\makeatother

\theoremstyle{plain}
\newtheorem{thm}{Theorem}[section]
\newtheorem{cor}[thm]{Corollary}
\newtheorem{lem}[thm]{Lemma}
\newtheorem{pro}[thm]{Proposition}

\theoremstyle{definition}

\numberwithin{equation}{section}

\newcounter{rmnum}

\def\symbolnote#1#2{\let\thefootn=\thefootnote%
\renewcommand{\thefootnote}{\fnsymbol{footnote}}%
\footnotemark[#1]%
\footnotetext[#1]{#2}%
\let\thefootnote=\thefootn
}

\newfont{\bbb}{msbm10 scaled \magstep1}
\newfont{\bbc}{msbm8 scaled \magstep0}

\newcommand{\N}{\mbox{\bbb N}}

\newcommand{\uin}{\mbox{\bbb I}}
\newcommand{\e}{\mbox{\rm e-dim}}

\begin{document}


\title{Extraordinary dimension of maps}
\author{A. Chigogidze}
\address{Department of Mathematical Sciences,
University of North Carolina at Greensboro,
P.O. Box 26170, Greensboro, NC 27402-6170, U.S.A.}
\email{chigogidze@uncg.edu}

\author{V. Valov}
\address{Department of Mathematics, Nipissing University,
100 College Drive, P.O. Box 5200, North Bay, ON, P1B 8L7, Canada}
\email{veskov@nipissingu.ca}
\thanks{The second named author was partially supported by his NSERC grant.}

 

\begin{abstract}
We establish a characterization of the extraordinary dimension of perfect maps between metrizable spaces. 
\end{abstract}

\maketitle

\markboth{A.~Chigogidze and V.~Valov}{Extraordinary dimension of maps}


\section{Introduction}

The paper deals with extensional dimension of maps, specially, with the extraordinary dimension introduced recently by \v{S}\v{c}epin \cite{es:98} and studied by the first author in \cite{ch:03}. If $L$ is  a $CW$-complex and $X$ a metrizable space, we write $\e X\leq L$ provided $L$ is an absolute extensor for $X$ (in such a case we say that the extensional dimension of $X$ is $\leq L$, see \cite{d:1}, \cite{d:2}). The extraordinary dimension of $X$ generated by a complex $L$, notation $\dim_LX$, is the smallest integer $n$ such that $\e X\leq\Sigma^nL$, where $\Sigma^nL$ is the $n$-th iterated suspension of $L$ (by $\Sigma^0L$ we always denote the complex $L$ itself). If $L$ is the $0$-dimensional sphere $S^0$, then $\dim_L$ coincides with the covering dimension $\dim$. We also write $\dim_Lf\leq n$, where $f\colon X\to Y$ is a given map, provided $\dim_Lf^{-1}(y)\leq n$ for every $y\in Y$. Next is our main result.

\begin{thm}\label{T:1.1}
Let $f\colon X\to Y$ be a $\sigma$-perfect map of metrizable spaces, let $L$ be a $CW$-complex and $n \geq 1$. Consider the following properties:

\begin{itemize}
\item[(1)] $\dim_Lf\leq n$;

\item[(2)] There exists an $F_{\sigma}$ subset $A$ of $X$ such that $\dim_L A\leq n-1$ and the restriction map $f|(X\backslash A)$ is of dimension $\dim_Lf|(X\backslash A) = 0$;

\item[(3)] There exists a dense and $G_{\delta}$ subset $\mathcal G$ of $C(X,\uin^n)$ with the source limitation topology such that $\dim_L(f\times g)= 0$ for every $g\in{\mathcal G}$;

\item[($3'$)] There exists a map  $g\colon X\to\uin^n$ is such that $\dim_L(f\times g) = 0$.
\end{itemize}

Then $(3)\Rightarrow (3')\Rightarrow (1)$ and  $(3')\Rightarrow (2)$. Moreover, $(1)\Rightarrow (3)$ provided $Y$ is a $C$-space and $L$ is countable.
\end{thm}

Here, $f\colon X\to Y$ is $\sigma$-perfect if $X$ is the union of countably many closed sets $X_i$ such that $f(X_i)\subset Y$ are closed and the restriction maps $f|X_i$ are perfect.

Theorem \ref{T:1.1} is inspired  by the following result of M. Levin and W. Lewis \cite[Theorem 1.8]{ll:02}: If $X$ and $Y$ are metrizable compacta then   $(3)\Rightarrow (3')\Rightarrow (1)$ and $(3)\Rightarrow (2')\Rightarrow (1)$, where $(2')$ is obtained from our condition $(2)$ by replacing $\dim_Lf|(X\backslash A)\leq 0$ with $\dim f|(X\backslash A)\leq 0$.
Moreover, the implication $(1)\Rightarrow (3)$ was also established in \cite{ll:02} for a finite-dimensional compactum $Y$  and a countable $CW$-complex $L$.

\smallskip
 
Therefore,  we have the following characterization of extraordinary dimension of perfect maps between metrizable spaces:

\begin{cor}\label{C:1.2}
Let $f\colon X\to Y$ be a perfect surjection between metrizable spaces with $Y$ being a $C$-space. If $L$ is a countable $CW$-complex, then the following conditions are equivalent:

\begin{itemize}
\item[(1)] $\dim_Lf\leq n$;

\item[(2)] There exists a dense and $G_{\delta}$ subset $\mathcal G$ of $C(X,\uin^n)$ with the source limitation topology such that $\dim_L(f\times g)\leq 0$ for every $g\in{\mathcal G}$;

\item[(3)] There exists a map  $g\colon X\to\uin^n$ is such that $\dim_L(f\times g)\leq 0$.
\end{itemize}

\smallskip\noindent
If, in addition, $X$ is compact, then each of the above three conditions is equivalent to the following one:

\begin{itemize}
\item[(4)] There exists an $F_{\sigma}$ set $A\subset X$ such that $\dim_L A\leq n-1$ and the restriction map $f|(X\backslash A)$ is of dimension $\dim f|(X\backslash A)\leq 0$.
\end{itemize}
\end{cor}

The equivalence of the first three conditions follow from Theorem \ref{T:1.1}. More precisely, by Theorem \ref{T:1.1} we have the following implications: 
$(2)\Rightarrow (3)\Rightarrow (1)\Rightarrow (2)$. When $X$ is compact, the result of Levin-Lewis which was mentioned above yields that 
$(2)\Rightarrow (4)\Rightarrow (1)$. Therefore, combining the last two chains of implications, we can obtain the compact version of Corollary \ref{C:1.2}.  

\medskip
Corollary \ref{C:1.2} is a parametric version of \cite[Theorem 4.9]{ch:03}. For the covering dimension $\dim$ such a characterization was obtained by Pasynkov \cite{bp:96} and Toru\'{n}czyk \cite{ht} in the realm of finite-dimensional compact metric spaces and extended in \cite{tv:02} to perfect maps between metrizable $C$-spaces. Since the class of $C$-spaces contains the class of finite-dimensional ones as a proper subclass (see \cite{re:95}),   
the compact version of Corollary \ref{C:1.2} is more general than the Levin-Lewis result  \cite[Theorem 1.8]{ll:02}.
It is interesting to know whether all the conditions $(1)$-$(4)$ in Corollary \ref{C:1.2} remain equivalent without the compactness requirement on $X$ and $Y$.

\medskip
The source limitation topology on $C(X,M)$, where $(M,d)$ is a metric space, can be described as follows:
a subset $U\subset C(X,M)$ is open if for every $g\in U$ there exists 
a continuous function $\alpha\colon X\to (0,\infty)$ such that $\overline{B}(g,\alpha)\subset U$. Here, $\overline{B}(g,\alpha)$ denotes the set  
$\{h\in C(X,M):d(g(x),h(x))\leq\alpha (x)\hbox{}~~\mbox{for each 
$x\in X$}\}$. 
The source limitation topology doesn't depend on the metric $d$ if $X$ is paracompact and $C(X,M)$ with this topology has the Baire property provided $(M,d)$ is a complete metric space.   Moreover, if $X$ is compact, then the source limitation topology coincides with the uniform convergence topology generated by $d$.  

All function spaces in this paper, if not explicitely stated otherwise, are equipped with the source limitation topology.


\section{Some preliminary results}

 Throughout this section $K$ is a closed and convex subset  of a given Banach space $E$ and   $f\colon X\to Y$  a  perfect map with $X$ and $Y$ paracompact spaces.
Suppose that for every $y\in Y$ we are given a property ${\mathcal P}(y)$ of maps $h\colon f^{-1}(y)\to K$ and let ${\mathcal P}=\{{\mathcal P}(y): y\in Y\}$.  By  
$C_{\mathcal P}(X|H,K)$ we denote the set of all bounded maps $g\colon X\to K$ such that $g|f^{-1}(y)$ has the property ${\mathcal P}(y)$ for every $y\in H$, where $H\subset Y$. We also consider the set-valued map $\displaystyle\psi_{\mathcal P}\colon Y\to 2^{C^*(X,K)}$, defined by the formula
$\psi_{\mathcal P}(y)=C^*(X,K)\backslash C_{\mathcal P}(X|\{y\},K)$, where $C^*(X,K)$ is the space of bounded maps from $X$ into $K$.

\begin{lem}
Suppose that $\mathcal P=\{{\mathcal P}(y)\}_{y\in Y}$ is a family of properties satisfying the following conditions:

\begin{itemize}
\item[(a)]  $C_{\mathcal P}(X|H,K)$ is open in $C^*(X,K)$ with respect to the source limitation topology for every closed $H\subset Y$;
\item[(b)] $g\in C_{\mathcal P}(X|\{y\},K)$ implies $g\in C_{\mathcal P}(X|U,K)$ for some neighborhood $U$ of $y$ in $Y$.
\end{itemize}
Then the map $\psi_{\mathcal P}$ has a closed graph provided $C^*(X,K)$ is equipped with the uniform convergence topology. 
\end{lem}

\begin{proof}
The proof of this lemma follows the arguments from the proof of \cite[Lemma 2.6]{tv:02}.
\end{proof}

Recall that a closed subset $F$ of the metrizable apace $M$ is said to be a $Z_m$-set in $M$,   if the set $C(\uin^m,M\backslash F)$ is dense in $C(\uin^m,M)$ with respect to the uniform convergence topology, where $\uin^m$ is the $m$-dimensional cube. If  
$F$ is a $Z_m$-set in $M$ for every $m\in\N$, we say that $F$ is a $Z$-set in $M$.

\begin{lem}
Suppose $y\in Y$ and ${\mathcal P}(y)$ satisfy the following condition:
\begin{itemize}
\item For every $m\in\N$  the set of all maps $h\in C(\uin^m\times f^{-1}(y),K)$ with each $h|(\{z\}\times f^{-1}(y))$, $z\in\uin^m$,  having the property $\mathcal P(y)$ $($as a map from $f^{-1}(y)$ into $K$$)$ is dense in $C(\uin^m\times f^{-1}(y),K)$ with respect to the uniform convergence topology.
\end{itemize}
Then, for every $\alpha\colon X\to (0,\infty)$ and $g\in C^*(X,K)$, $\psi_{\mathcal P}(y)\cap\overline{B}(g,\alpha)$  is a $Z$-set in 
$\overline{B}(g,\alpha)$ provided $\overline{B}(g,\alpha)$ is considered as subset of $C^*(X,K)$ equipped with the uniform convergence topology and
$\psi_{\mathcal P}(y)\subset C^*(X,K)$ is closed.
\end{lem}

\begin{proof}
See the proof of \cite[Lemma 2.8]{tv:02}
\end{proof}

\begin{pro}
Let  $Y$ be a $C$-space and  $\mathcal P=\{\mathcal P(y)\}_{ y\in Y}$ such that:

\begin{itemize}
\item[(a)] the map $\psi_{\mathcal P}$ has a closed graph;
\item[(b)] $\psi_{\mathcal P}(y)\cap\overline{B}(g,\alpha)$ is a $Z$-set in $\overline{B}(g,\alpha)$
for any continuous function $\alpha\colon X\to (0,\infty)$, $y\in Y$ and $g\in C^*(X,K)$,  where 
$\overline{B}(g,\alpha)$ is considered as a subspace of $C^*(X,K)$ with the uniform convergence topology. 
\end{itemize}
Then the set $\{g\in C^*(X,K): g\in C_{\mathcal P}(X|\{y\},K)\hbox{}~\mbox{for every $y\in Y$}\}$ is dense in $C^*(X,K)$ with respect to the source limitation topology.
\end{pro}

\begin{proof}
Let $G=\{g\in C^*(X,K): g\in C_{\mathcal P}(X|\{y\},K)\hbox{}~\mbox{for every $y\in Y$}\}$.
It suffices to show that, for fixed $g_0\in C^*(X,K)$ and a positive continuous function $\alpha\colon X\to (0,\infty)$, there exists $g\in \overline{B}(g_0,\alpha)\cap G$. We equip $C^*(X,K)$ with the uniform convergence topology and 
consider
the constant (and hence, lower semi-continuous) convex-valued map $\phi\colon Y\to 2^{C^*(X,K)}$, 
$\phi(y)=\overline{B}(g_0,\alpha_1)$, where $\alpha_1(x)=\min\{\alpha (x), 1\}$.
Because of the conditions (a) and (b), we  
can apply the selection theorem \cite[Theorem 1.1]{gv:99} to obtain a  continuous map $h\colon Y\to C^*(X,K)$ such that  $h(y)\in\phi(y)\backslash\psi_{\mathcal P}(y)$ for every $y\in Y$.
Observe that $h$ is a map from $Y$ into $\overline{B}(g_0,\alpha_1)$ such that $h(y)\in C_{\mathcal P}(X|\{y\},K)$ for every $y\in Y$. Then
$g(x)=h(f(x))(x)$, $x\in X$, defines a bounded map $g\in \overline{B}(g_0,\alpha)$ such that  $g|f^{-1}(y)=h(y)|f^{-1}(y)$, $y\in Y$. Therefore, $g\in C_{\mathcal P}(X|\{y\},K)$ for all $y\in Y$, i.e., $g\in \overline{B}(g_0,\alpha)\cap G$.
\end{proof}

\section{Proof of Theorem \ref{T:1.1}}

\noindent
$\bf{(1)\Rightarrow (3)}$ Suppose that $L$ is countable and $Y$ is a $C$-space.  Let $X_i$ be closed subsets of $X$ such that each $f_i=f|X_i\colon X_i\to Y_i=f(X_i)$ is a perfect map and $Y_i$ is closed in $Y$. Then all $Y_i$'s are $C$-spaces, and since the restriction maps $\pi_i\colon C(X,\uin^n)\to C(X_i,\uin^n)$, $\pi_i(g)=g|X_i$, are open, the proof of this implication is reduced to the case when $f$ is a perfect map. Consequently, we may assume that $f$ is perfect.
  
By \cite[Theorem 1.1]{tv1} (see also \cite{bp:98}), there exists a map $q$ from $X$ into the Hilbert cube $Q$ such that $f\times q\colon X\to Y\times Q$ is an embedding. Let $\{W_i\}_{i\in\N}$ be  a countable finitely-additive base for $Q$. 
For every $i$ we choose a sequence of mappings $h_{ij}\colon\overline{W_i}\to L$, representing all the homotopy classes of mappings from $\overline{W_i}$ to $L$ (this is possible because $L$ is a countable $CW$-complex).  Following the notations from Section 2, for fixed $i$, $j$ and $y\in Y$ we say that a map $g\in C(X,\uin^n)$ has the property ${\mathcal P}_{ij}(y)$ provided

\smallskip\noindent
 the map $h_{ij}\circ q\colon q^{-1}(\overline{W_{i}})\to L$ can be continuously extended to a map over the set
$q^{-1}(\overline{W_{i}})\cup\big( f^{-1}(y)\cap g^{-1}(t)\big)$ for every $t\in g(f^{-1}(y))$.

\smallskip\noindent 
Let ${\mathcal P}_{ij}=\{{\mathcal P}_{ij}(y): y\in Y\}$ and for every 
$H\subset Y$ we denote $\displaystyle C_{{\mathcal P}_{ij}}(X|H,\uin^n)$ by $C_{ij}(X|H,\uin^n)$. Hence, $C_{ij}(X|H,\uin^n)$ consists of all $g\in C(X,\uin^n)$ having the property $C_{ij}(y)$ for every $y\in H$. Let $\psi_{ij}\colon Y\to 2^{C(X,\uin^n)}$ be the set-valued map
$\psi_{ij}(y)=C(X,\uin^n)\backslash C_{ij}(X|\{y\},\uin^n)$.

\begin{lem}
Let $g\in C_{ij}(X|\{y\},\uin^n)$. Then, there exist a neighborhood $U_y$ of $y$ in $Y$ and a neighbourhood $V_t\subset\uin^n$ of each 
$t\in g(f^{-1}(y))$ such that  $h_{ij}\circ q$ can be extended to a map from $q^{-1}(\overline{W_{i}})\cup\big( f^{-1}(U_y)\cap g^{-1}(V_t)\big)$ into $L$.
\end{lem}

\begin{proof}
Since $g\in C_{ij}(X|\{y\},\uin^n)$, $h_{ij}\circ q$ can be extended to a map from $q^{-1}(\overline{W_{i}})\cup\big( f^{-1}(y)\cap g^{-1}(t)\big)$ into $L$ for every $t\in g(f^{-1}(y))$.  Because $L$ is an absolute neighborhood extensor for $X$, there exists and open set $G_t\subset X$ containing $f^{-1}(y)\cap g^{-1}(t)$ and a map $h_t\colon q^{-1}(\overline{W_{i}})\cup G_t\to L$ extending  $h_{ij}\circ q$. Using that $f\times g$ is a closed map, we can find a neighborhood $U_y^t\times V_t$ of $(y,t)$ in $Y\times\uin^n$ such that $S_t=(f\times g)^{-1}(U_y^t\times V_t)\subset G_t$.  Next, choose finitely many points  $t(k)$, $k=1,2,.,m$, with 
$f^{-1}(y)\subset\bigcup_{k=1}^{k=m} S_{t(k)}$ and a neighbothood $U_y$ of $y$ in $Y$ such that $U_y\subset\bigcap_{k=1}^{k=m}U_y^{t(k)}$ and $f^{-1}(U_y)\subset \bigcup_{k=1}^{k=m} S_{t(k)}$ (this can be done since $f$ is perfect).  
 If $t\in g(f^{-1}(y))$, then $t\in V_{t(k)}$ for some $k$ and $f^{-1}(U_y)\cap g^{-1}(V_{t(k)})\subset S_{t(k)}$.  Since , $S_{t(k)}\subset G_{t(k)}$, the map $h_{t(k)}$ is an extension of $h_{ij}\circ q$ over the set $q^{-1}(\overline{W_{i}})\cup\big( f^{-1}(U_y)\cap g^{-1}(V_{t(k)})\big)$
\end{proof}

\begin{lem}
The set $C_{ij}(X|H,\uin^n)$ is open in $C(X,\uin^n)$ for any $i,j$ and closed $H\subset Y$.
\end{lem}

\begin{proof}
We follow the proof of \cite[Lemma 2.5]{tv:02}. For a fixed $g_0\in C_{ij}(X|H,\uin^n)$ we are going to find a function $\alpha\colon X\to (0,\infty)$ such that $\overline{B}(g_0,\alpha)\subset C_{ij}(X|H,\uin^n)$. By Lemma 3.1, for every $z=(y,t)\in (f\times g_0)((f^{-1}(H))$ there exists a neighborhood $U_z$ in $Y\times\uin^n$ such that 

\medskip\noindent
(1)\hspace{0.4cm}$h_{ij}\circ q$ can be extended to a map from $q^{-1}(\overline{W_{i}})\cup (f\times g_0)^{-1}(U_z)$ into $L$.  

\medskip\noindent
Obviously, $K=(f\times g_0)((f^{-1}(H))$ is closed in $Y\times\uin^n$, so there exists open $G\subset Y\times\uin^n$ with $K\subset G\subset\overline{G}\subset U=\bigcup\{U_z: z\in K\}$. Then $\nu=\{U_z: z\in K\}\cup\{(Y\times\uin^n)\backslash\overline{G}\}$ is an open cover of $Y\times\uin^n$. Let $\gamma$ be an open locally finite cover of $Y\times\uin^n$ such that the family 

\medskip\noindent
(2)\hspace{0.4cm}$\{St(W,\gamma): W\in\gamma\}$ refines $\nu$ and $St(W,\gamma)\subset G$ provided $W\cap K\neq\emptyset$. 

\medskip\noindent
Consider the metric $\rho=d + d_1$ on $Y\times\uin^n$, where $d$ is a metric on $Y$ and $d_1$  the usual metric on $\uin^n$, and define the function $\alpha\colon X\to (0,\infty)$ by $\alpha (x)=2^{-1}\sup\{\rho\big((f\times g_0)(x), (Y\times\uin^n)\backslash W\big): W\in\gamma\}$.  Let show that $\overline{B}(g_0,\alpha)\subset C_{ij}(X|H,\uin^n)$.   
Take $g\in \overline{B}(g_0,\alpha)$,  
$y\in H$ and $t\in g(f^{-1}(y))$. Then, $(y,t)=(f\times g)(x)$ for some $x\in f^{-1}(y)$. Since $g$ is $\alpha$-close to $g_0$, 
there exists $W\in\gamma$ such that $W\cap K\neq\emptyset$ and $W$ contains both $(f\times g)(x)$ and $(f\times g_0)(x)$. 
It follows from (2) that $(f\times g)^{-1}(W)\subset (f\times g_0)^{-1}(U_z)$ for some $z\in K$.  In particular, 
$f^{-1}(y)\cap g^{-1}(t)\subset (f\times g_0)^{-1}(U_z)$.  Consequently, by (1), $h_{ij}\circ q$ is extendable to a map from 
$q^{-1}(\overline{W_{i}})\cup\big(f^{-1}(y)\cap g^{-1}(t)\big)$ into $L$.  Therefore,  $\overline{B}(g_0,\alpha)\subset C_{ij}(X|\{y\},\uin^n)$ for every $y\in H$ which completes the proof.
\end{proof}

Because of Lemma 3.1 and Lemma 3.2, we can apply Lemma 2.1 to obtain the following corollary.

\begin{cor}
For any $i$ and $j$ the map $\psi_{ij}$ has a closed graph.
\end{cor} 

\begin{lem}
Let $g\in C(X,\uin^n)$, $\alpha\colon X\to (0,\infty)$ and $y\in Y$. Then, for any $i,j$, $\psi_{ij}(y)\cap\overline{B}(g,\alpha)$  is a $Z$-set in 
$\overline{B}(g,\alpha)$ provided $\overline{B}(g,\alpha)$ is considered as a subset of $C(X,\uin^n)$ with the uniform convergence topology.
\end{lem}

\begin{proof}
It follows from  \cite[Theorem 1.8, $(1)\Rightarrow (3)$]{ll:02} that 
if  $m\in\N$, then all maps $g\colon\uin^m\times f^{-1}(y)\to\uin^n$ such that $\e\big((\{z\}\times f^{-1}(y))\cap g^{-1}(t)\big)\leq L$ for every $z\in\uin^m$ and $t\in\uin^n$, form a dense subset $G$ of $C(\uin^m\times f^{-1}(y))$ with the uniform convergence topology. 
It is clear that, for every $g\in G$ and $z\in\uin^m$, the restriction $g|\big(\{z\}\times f^{-1}(y)\big)$, considered as a map on $f^{-1}(y)$, has the following property: $h_{ij}\circ q$ can be extended to a map from $q^{-1}(\overline{W_{i}})\cup (f^{-1}(y)\cap g^{-1}(t))$ into $L$ for any $t\in\uin^n$.  Hence, we can apply Lemma 2.2 to conclude that $\psi_{ij}(y)\cap\overline{B}(g,\alpha)$  is a $Z$-set in 
$\overline{B}(g,\alpha)$.
\end{proof}

Now, we can finish the proof of this implication. Because of Corollary 3.3 and Lemma 3.4, we can apply Proposition 2.3 to obtain that the set $C_{ij}=C_{ij}(X|Y,\uin^n)$ is dense in $C(X,\uin^n)$ for every $i,j$. Since, by Lemma 3.2, all $C_{ij}$ are also open, their intersection $\mathcal G$ is dense and $G_{\delta}$ in $C(X,\uin^n)$.  Let show that $\dim_L (f\times g)\leq 0$
for every $g\in\mathcal G$, i.e., $\e(f\times g)\leq L$. We fix $y\in Y$ and $t\in\uin^n$ and consider the fiber $(f\times g)^{-1}(y,t)=f^{-1}(y)\cap g^{-1}(t)$.  Take a closed set $A\subset f^{-1}(y)\cap g^{-1}(t)$ and a map $h\colon A\to L$. Because the map $q_y=q|f^{-1}(y)$ is a homeomorphism, $h^{'}=h\circ q_y^{-1}\colon q(A)\to L$ is well defined. Next, extend $h^{'}$ to a map from a neighborhood  $W$ of $q(A)$ (in $Q$) into $L$ and find $W_i$ with $q(A)\subset W_i\subset\overline{W_i}\subset W$. Therefore, there exists a map
$h^{''}\colon\overline{W_i}\to L$ extending $h^{'}$.  Then $h^{''}$ is homotopy equivalent to some $h_{ij}$, so are 
$h^{''}\circ q$ and $h_{ij}\circ q$ (considered as maps from $q^{-1}(\overline{W_{i}})$ into $L$).  Since $h_{ij}\circ q$ can be extended to a map from $q^{-1}(\overline{W_{i}})\cup\big( f^{-1}(y)\cap g^{-1}(t)\big)$ into $L$, by the Homotopy Extension Theorem, there exists a map $\bar{h}\colon q^{-1}(\overline{W_{i}})\cup\big( f^{-1}(y)\cap g^{-1}(t)\big)\to L$ extending $h^{''}\circ q$. Obviously, $\bar{h}|\big( f^{-1}(y)\cap g^{-1}(t)\big)$ extends  $h$. Hence, 
$\e\big( f^{-1}(y)\cap g^{-1}(t)\big)\leq L$. 

\bigskip\noindent
$\bf{(3)\Rightarrow (3')\Rightarrow (1)}$ The implication $(3)\Rightarrow (3')$ is trivial.  It is easily seen that in the proof of $(3')\Rightarrow (1)$ we can assume $f$ is perfect. Let 
$g\colon X\to\uin^n$ be such that $\dim_L(f\times g)\leq 0$ and $y\in Y$. Then the restriction $g|f^{-1}(y)\colon f^{-1}(y)\to\uin^n$ is a pefect map with all of its fibers having extensional dimension $\e\leq L$. Hence, by \cite[Corollary]{cv:1}, $\e f^{-1}(y)\leq\Sigma^nL$, i.e, $\dim_Lf\leq n$. 

\bigskip\noindent
$\bf{(3')\Rightarrow (2)}$ Because of the countable sum theorem, we can suppose that $f$ is perfect. We fix a map  
$g\colon X\to\uin^n$ such that $\dim_L(f\times g)\leq 0$. According to \cite[Lemma 4.1]{tv:02}, there exists an $F_{\sigma}$ subset $B\subset Y\times\uin^n$ such that $\dim B\leq n-1$ and $\dim (\{y\}\times\uin^n)\backslash B\leq 0$ for every $y\in Y$. Then, applying again \cite[Corollary]{cv:1}, we conclude that the set $A=(f\times g)^{-1}(B)$ is as required.  
   

\bigskip

\end{document}